\documentclass[11pt]{article}
\usepackage{amsmath, amsthm, amssymb}    
\usepackage{bridges2}                     
\usepackage{graphicx}                    

\usepackage[hyphens]{url}
\usepackage[colorlinks=true, urlcolor=blue, citecolor=black, linkcolor=black]{hyperref}                     

\usepackage{enumitem}
\usepackage{color}
\usepackage{wrapfig}
\usepackage{pdflscape}
\usepackage{array}
\usepackage{MnSymbol}
\usepackage{pinlabel}
\usepackage[font=small]{caption, subcaption}
\usepackage{floatrow}

\usepackage{lipsum}

\newcommand{\ZZ}{\mathbb{Z}}
\newcommand{\RR}{\mathbb{R}}

\renewcommand{\th}{^{\textrm{th}}}

\title{On the Group Structure of ``Magic'' Leatherworking Braids}

\author{William Hobkirk\textsuperscript{1}, Abigail Hollingsworth\textsuperscript{2}, Elisabetta A. Matsumoto\textsuperscript{3}, and Corbin Reid\textsuperscript{4}\\ 
\vspace{-2pt}\\
\textsuperscript{1}School of Mathematics and Statistics, University of Sydney, Sydney, Australia;\\ { william.hobkirk@sydney.edu.au}\\
\textsuperscript{2}Warwick Mathematics Institute, University of Warwick, Coventry, UK; \\{ abigail.hollingsworth@warwick.ac.uk}\\
\textsuperscript{3}School of Physics, Georgia Institute of Technology, Atlanta, GA, USA;\\ { sabetta@gatech.edu}\\
\textsuperscript{4}School of Mathematics, Monash University, Melbourne, Australia; \\{corbin.reid1@monash.edu}
\vspace{-13pt}
}
\date{}				

\begin{document}

\maketitle

\thispagestyle{empty}

\begin{abstract}
\emph{Magic braids} are used in leatherworking to make intricate straps and bands. They appear to be impossible to make. To determine which braids can be made with this leatherworking technique, we explore their relation to several braid groups.
\end{abstract}


\vspace{-12pt}
\section*{Introduction}

\begin{figure}[htb]
\centering
\includegraphics[height=2.5in]{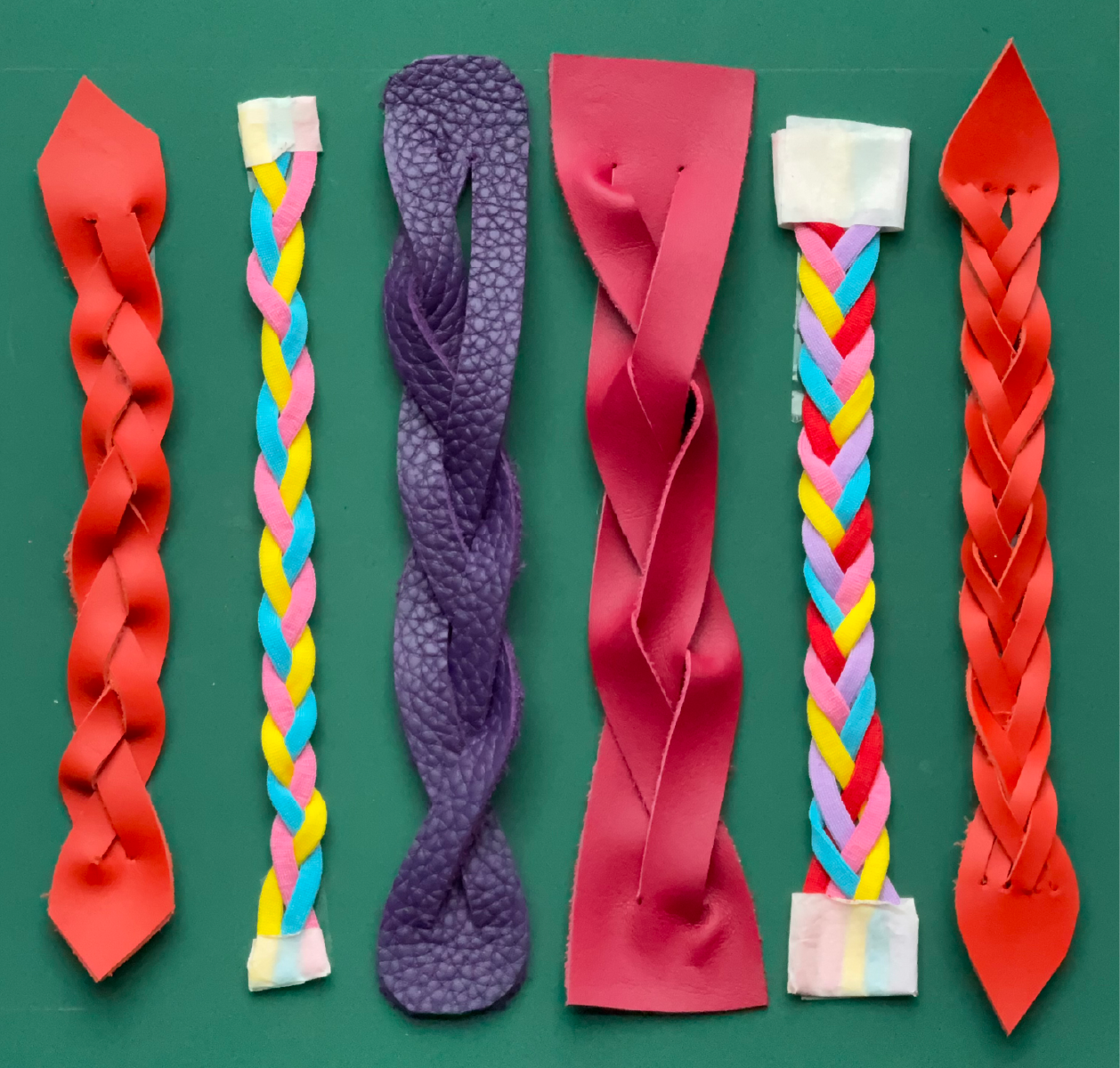}
\caption{{An example of several different} \emph{magic braids}. }
\label{Fig:examples}
\end{figure}

In leatherworking, \emph{magic braids}, also known as \emph{mystery braids}, are crafted by starting with a rectangular piece of leather with $n-1$ slits cut to form $n$ strands. This object is then manipulated by the leatherworker without cutting, tearing or stretching the leather and in this way they are able to move the strips of leather into many different configurations. In most of these configurations some of the strips are twisted. The magic lies in the leatherworker's ability to work the strands into a braid that is non-trivial but nevertheless lies flat without any twists in the strands. Figure~\ref{Fig:examples} shows several examples of magic braids. These braids appear to be \emph{magic} because they seem impossible to make until one reduces the situation to understanding how the twists and crossings of the leather strips change as we apply certain moves.

In topology, embeddings of an object in space are often considered equivalent if there exists continuous distortion of the embedding space that carries one embedding onto the other. This is called an \emph{ambient isotopy}. The relevant embedding space here is the three-sphere $S^3$ and the embedded object is the cut rectangle of leather. Different configurations of this cut rectangle are different embeddings. The manipulations of the leatherworker are ambient isotopies of this space since the leatherworker performs no destructive actions like cutting or pasting pieces of leather together. In this sense the configurations that a leatherworker is able to achieve are limited to those that are ambient isotopic to the initial cut rectangle. The one caveat to this is that since leather is not stretchy like the usual topological analogies of rubber or plasticene, we require that the order of the strips meeting at the \emph{coupons} at the top and bottom of the leather braid be fixed by these ambient isotopies. The goal of our mathematical study of magic braids is to clear up the relationship between these leather objects and the well understood \emph{braid groups} so that we can determine which leather braids are possible to make. Magic braids have been previously discussed by Murasugi and Kurpita in 1999 \cite{murasugi1999}, where they were referred to as \emph{Mexican plaits}, and Diamond 2025 \cite{diamond2025thesis}, where they are referred to as \emph{mystery braids}. In this paper, we introduce \emph{framing} to the magic braid formalism, where each strand in the braid is taken to be a ribbon, rather than a simple curve. While the framing can make the algebra more challenging, we change the setting from the traditional \emph{planar} braid group to the braid group in a thickened sphere to make the algebra tractable.

\begin{figure}[htb]
\centering
\includegraphics[height=2in]{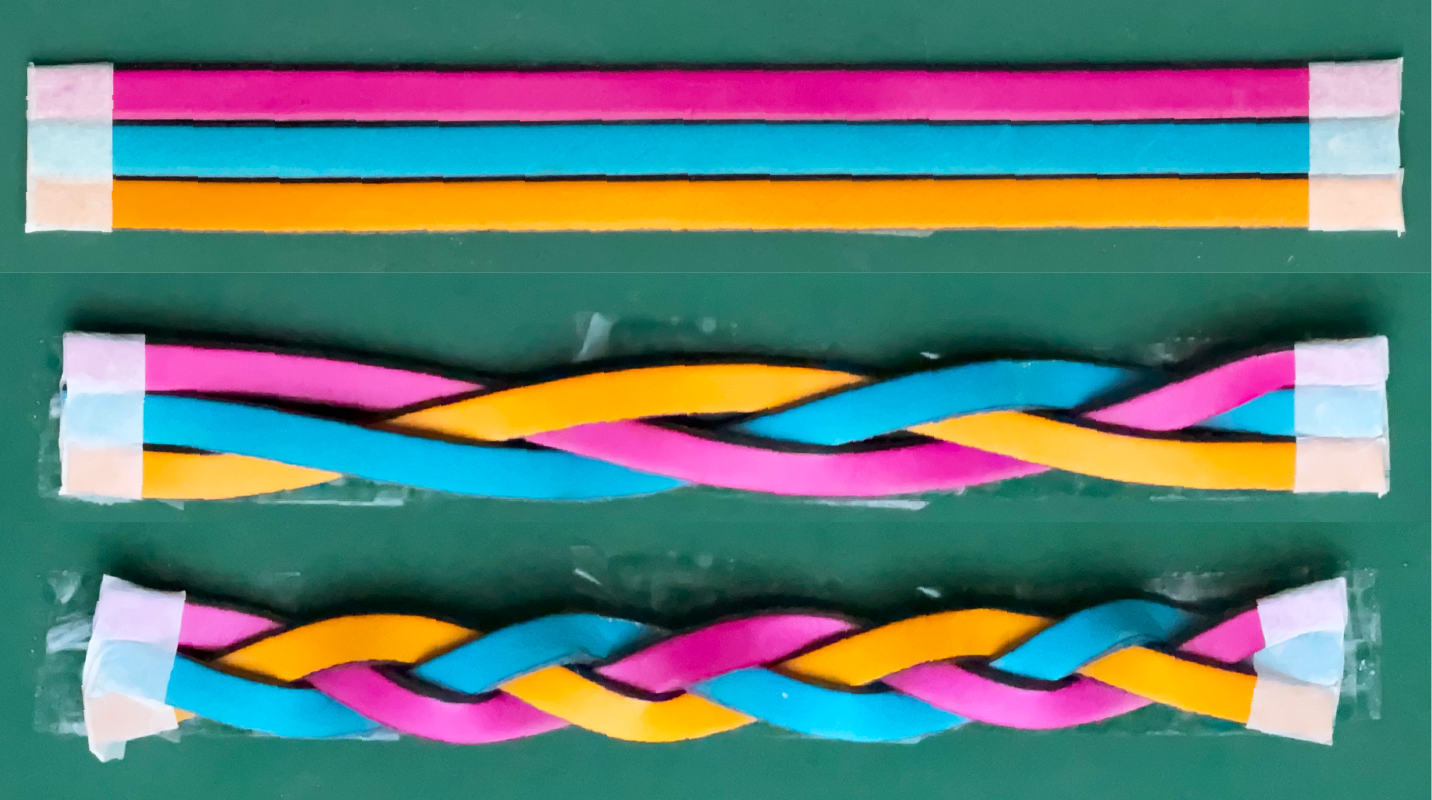}
\caption{{There exists an ambient isotopy from a leather strip to any magic braid. }}
\label{Fig:isotopy}
\end{figure}

\vspace{-18pt}
\section*{Planar Braid Groups: Their Presentations and Relations}

Before we get to magic braids, we start by briefly describing the \emph{braid group}, the \emph{pure braid group}, the \emph{framed braid group} and the \emph{pure framed braid group}. For more thorough discussions of the braid group and the pure braid group, see \cite{birmanbraids}, and for the framed braid group, see \cite{ko1992framedbraids}.

\begin{wrapfigure}[14]{l}{1.89in}
\vspace{-12pt}
    \begin{minipage}{.85in}
    \centering\captionsetup[subfigure]{justification=centering}
     \includegraphics[width=.85in]{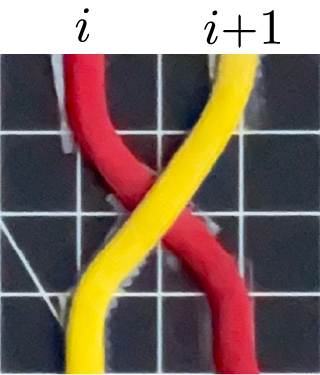}
    \vspace{-17pt}
    \subcaption{$\sigma_{i}$}
    \label{fig:Si}
    \vspace{1pt}

    \includegraphics[width=.85in]{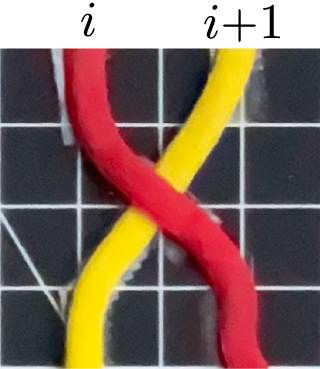}
    \vspace{-16pt}
    \subcaption{$\sigma_{i}^{-1}$}
    \label{fig:Si_inv}
   \end{minipage}
   \
   \begin{minipage}{.85in}
    \centering\captionsetup[subfigure]{justification=centering}
    \includegraphics[width=.85in]{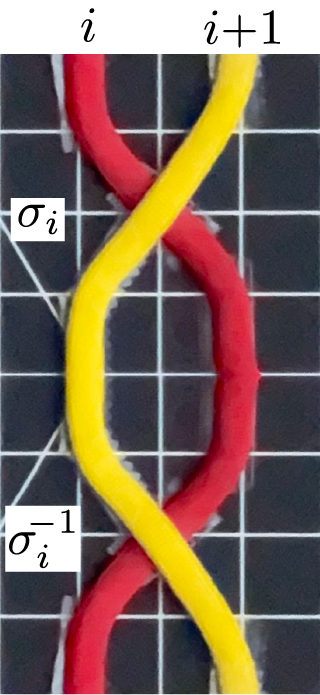}
    \vspace{-16pt}
    \subcaption{$\sigma_{i}\sigma_{i}^{-1}=\textrm{id}$}
    \label{fig:Si_Si_inv}
   \end{minipage}

\vspace{-8pt}
\caption{{Braid group generators.}}\label{fig:bg_gen}
\end{wrapfigure}

\noindent \textbf{\emph{Braid Group}}

\begin{figure}
\centering
\begin{minipage}{3.6in}
\centering\captionsetup[subfigure]{justification=centering}
\includegraphics[height=1.2in]{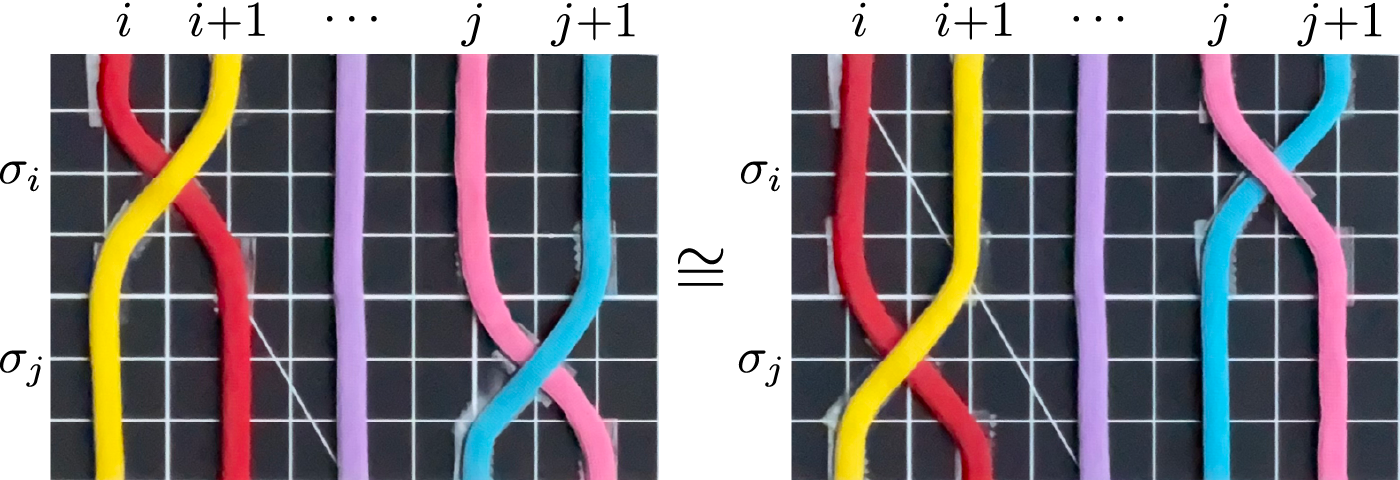}
\subcaption{$\sigma_{i}\sigma_{j}  = \sigma_{j}\sigma_{i},$ \emph{for} $| i-j|>1$}
\label{fig:bg_r1}
\end{minipage}
\begin{minipage}{2.83in}
\centering\captionsetup[subfigure]{justification=centering}
\includegraphics[height=1.5in]{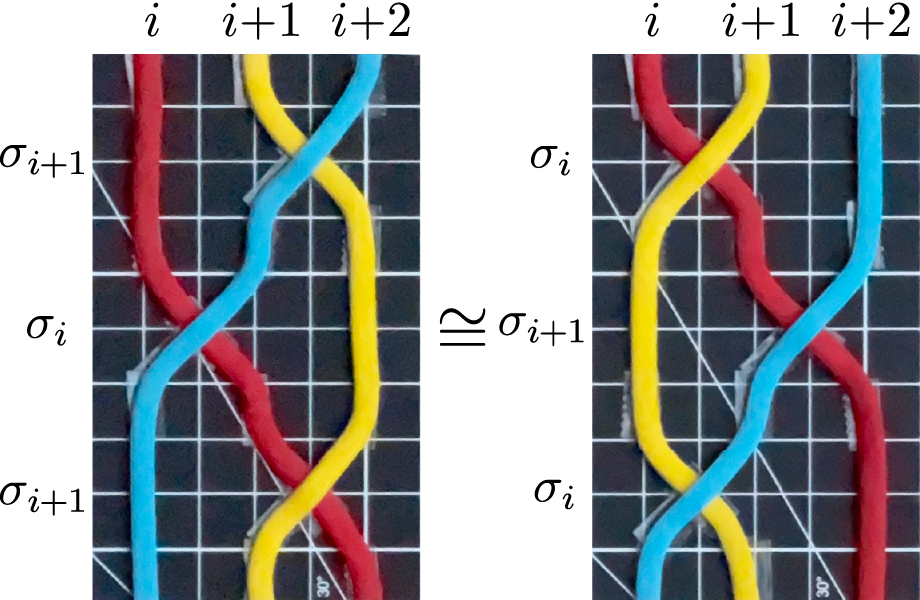}
\vspace{-3pt}
\subcaption{$\sigma_{i}\sigma_{i+1}\sigma_{i} = \sigma_{i+1}\sigma_{i}\sigma_{i+1}$}
\label{fig:bg_r2}
\end{minipage}

\vspace{5pt}
\caption{{Braid group relations.}}
\label{fig:bg_gen}
\vspace{-15pt}
\end{figure}

\noindent 
The standard \emph{braid group} on $n$ strands, consisting of braids in the thickened plane $\RR^2\times [0,1]$, will be referred to here as the \emph{planar braid group} $B_n^{\mathbb{R}^2}$ to distinguish from variations used later. 
A standard presentation of this group --- the Artin presentation~\cite{artin1947braids} --- uses the generating set $\{\sigma_i\}_{i=1}^{n-1}$, where $\sigma_i$ is a crossing of strand $i+1$ over strand $i$, as in Figure~\ref{fig:Si}. The inverse $\sigma_i^{-1}$ instead has strand $i$ crossing over strand $i+1$, as shown in Figure~\ref{fig:Si_inv}. These are both examples of a \emph{braid diagram}, in which the braid is projected to a vertical plane in $\RR^2\times[0,1]$. Figure~\ref{fig:Si_Si_inv} shows the composition $\sigma_i \sigma_i^{-1}$, which is equivalent to the identity braid as the crossings can be removed by isotopy (a Reidemeister II move). 
Figures~\ref{fig:bg_r1} (the crossing slide) and~\ref{fig:bg_r2} (the Reidemeister III relation) visualise the complete set of relations in $B_n$. 

\noindent\textbf{\emph{Pure Braid Group}}

\noindent Though, in general, the ordering of the strands at the top and bottom of a braid may be some nontrivial permutation of one another, both ends of a magic braid are fixed and keep the same ordering. We can construct a homomorphism from the braid group $B^{\RR^2}_n$ to the \emph{symmetric group} on $n$ elements $S_n$ by observing that $\sigma_i$ transposes strands $i$ and $i+1$, and that $S_n$ is generated by such transpositions. The kernel of the homomorphism that sends each $\sigma_i$ to the element of $S_n$ that interchanges $i$ and $i+1$ consists of braids where the strands at either end are not permuted with respect to one another. This kernel is a subgroup  $P_n$ of the braid group on $n$ strands, called the \emph{pure braid group}.

\begin{figure}[h!]
\centering
\captionsetup[subfigure]{justification=centering}
\begin{minipage}{1.75in}
\includegraphics[width=\linewidth]{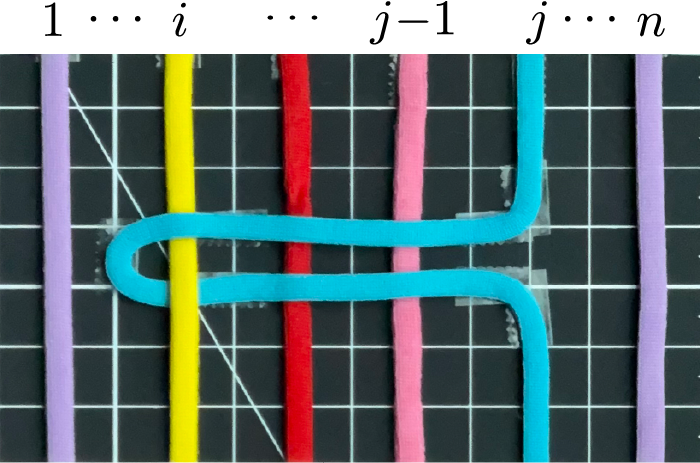}
\subcaption{$A_{ij}$}
\label{fig:pbg_Aij}
\end{minipage}
\ \
\begin{minipage}{1.75in}
\includegraphics[width=\linewidth]{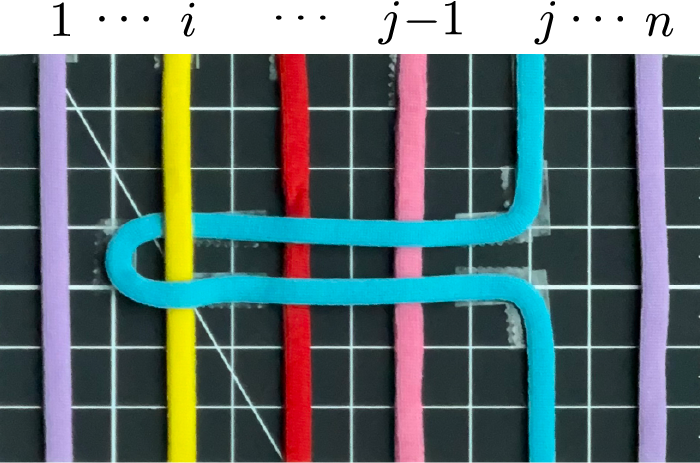}
\subcaption{$A_{ij}^{-1}$}
\label{fig:pbg_Aij_inv}
\end{minipage}
\ \
\begin{minipage}{1.75in}
\includegraphics[width=\linewidth]{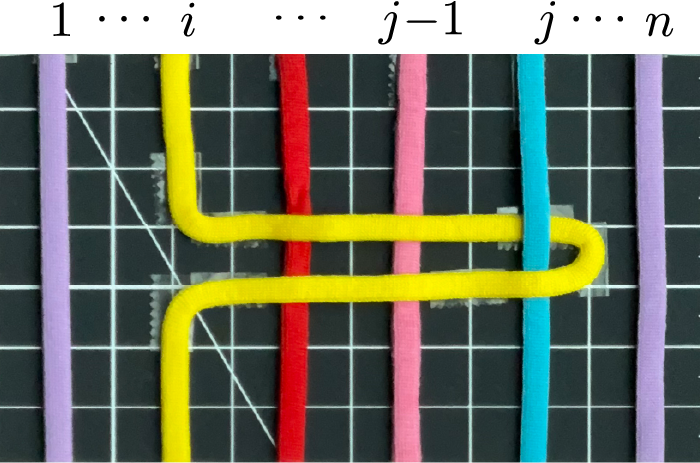}
\subcaption{$A_{ij} = A_{ji}^{-1}$} 
\label{fig:pbg_Aji}
\end{minipage}

\caption{{Generators for the pure braid group} $P_n$.}
\label{Fig:pbg_gen}
\vspace{-5pt}
\end{figure}

There is a presentation of $P_n$ with generators $A_{ij}$, which each pass the $j\th$ strand over the $i\th$ strand and back underneath to the $j\th$ position; see Figure~\ref{fig:pbg_Aij}. The inverse $A_{ij}^{-1}$ passes the $j\th$ strand first under, then over the $i\th$ strand as in Figure~\ref{fig:pbg_Aij_inv}. This generator can be expressed in terms of the generators of $B_n$ as $A_{ij}= \sigma_{j-1}\sigma_{j-2} \dots \sigma_{i+1}\sigma_i^2\sigma_{i+1}^{-1}\dots\sigma_{j-2}^{-1}\sigma_{j-1}^{-1},$ where $1\le i<j \le n$. 
From Figure~\ref{fig:pbg_Aji}, we can see that there is an isotopy that sends $A_{ji}$ to $A_{ij}$ for the case that $i>j$. 

This presentation of the pure braid group on $n$ strands $P_n$ has relations:
\begin{equation}
A_{rs}^{-1}A_{ij}A_{rs} = \left\{
\begin{array}{ll}
A_{ij} & \textrm{if $s<i$ or $i<r<s<j$}\\
A_{rj}A_{ij}A_{rj}^{-1} & \textrm{if $s=i$}\\
A_{rj}A_{sj}A_{ij}(A_{rj}A_{sj})^{-1} & \textrm{if $i=r<s<j$}\\
A_{rj}A_{sj}A_{rj}^{-1} A_{sj}^{-1} A_{ij} (A_{rj}A_{sj}A_{rj}^{-1} A_{sj}^{-1})^{-1} & \textrm{if $r<i<s<j$}
\end{array} 
\right.
\label{eq:pbg_rel}
\end{equation}

\subsection*{Framed Braid Groups}

\begin{wrapfigure}[7]{r}{1.75in}
\vspace{-37pt}
\centering
    \centering\captionsetup[subfigure]{justification=centering}
    \includegraphics[width=1.75in]{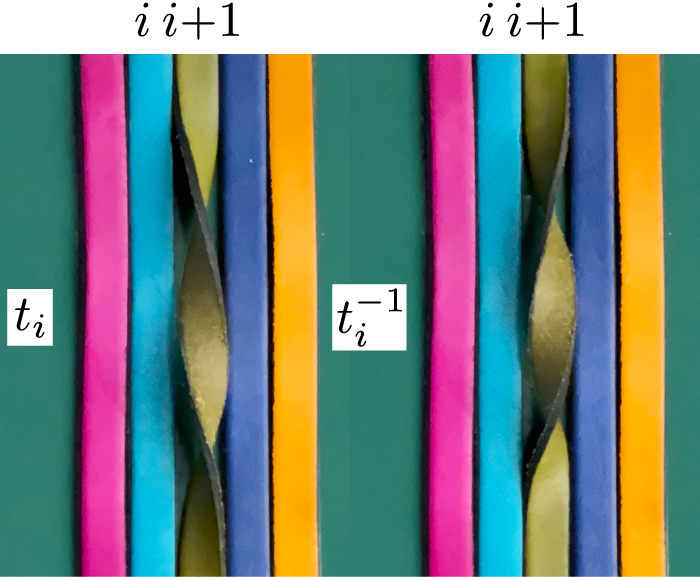}
    \vspace{-18pt}
 
\vspace{-5pt}
\caption{{Twists:} $t_i$ \emph{and} $t_i^{-1}$.}
\label{fig:fbg_gen}
\end{wrapfigure}
One of the reasons magic braids look so striking, particularly when they are made from leather, is that the shiny outside of the leather is always facing the same direction along the braid, while the suede/grain side is always facing the opposite side of the braid. When making a magic braid, the leatherworker often has to take care to smooth out the strips in the braid to ensure this appearance. We can think of each of the strips in the identity magic braid as a ribbon that can also be twisted. Because the top of the braid and the bottom of the braid must have the same orientation, only full twists on a strand can be consistent with the identity magic braid. However, in manipulating the strands, we may introduce unwanted twist, which will need to be canceled out by later moves.

\begin{figure}[t!]
\vspace{-0pt}
    \begin{minipage}{2in}
    \centering\captionsetup[subfigure]{justification=centering}
     \includegraphics[width=1.75in]{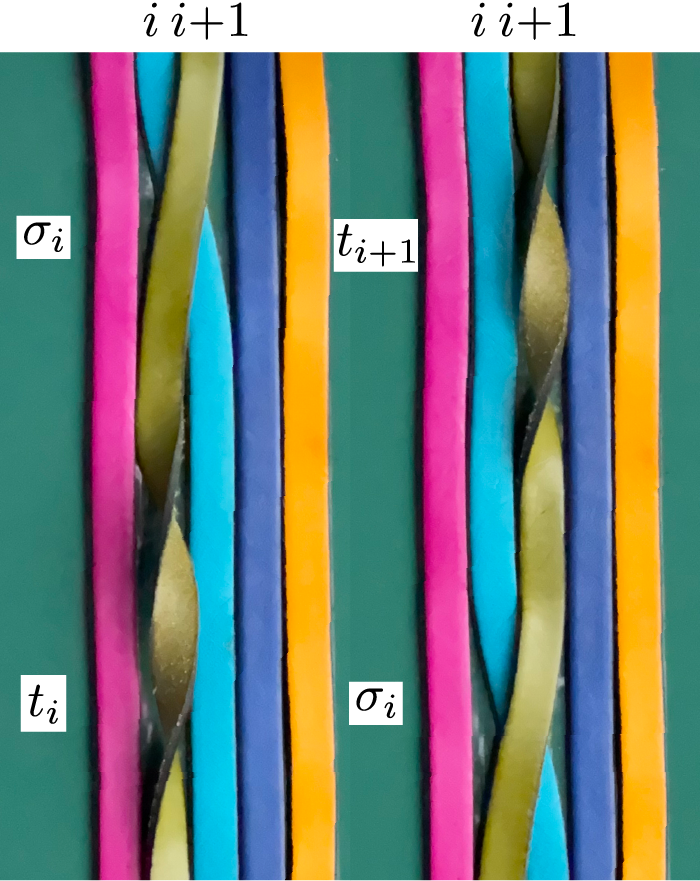}
    \subcaption{$\sigma_i t_i = t_{i+1}\sigma_i$}
    \label{fig:fbg_r1}
    \vspace{3pt}    
   \end{minipage}
   \
   \begin{minipage}{2in}
    \centering\captionsetup[subfigure]{justification=centering}
    \includegraphics[width=1.75in]{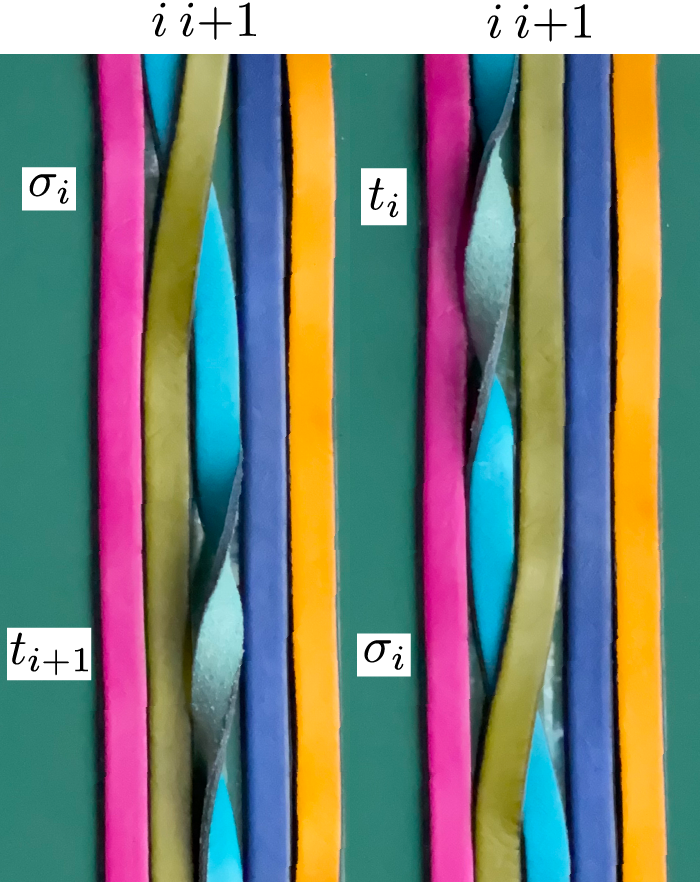}
    \subcaption{$\sigma_i t_{i+1} = t_{i}\sigma_i$}
    \label{fig:fbg_r2}
   \end{minipage}
     \
   \begin{minipage}{2in}
    \centering\captionsetup[subfigure]{justification=centering}
    \includegraphics[width=1.75in]{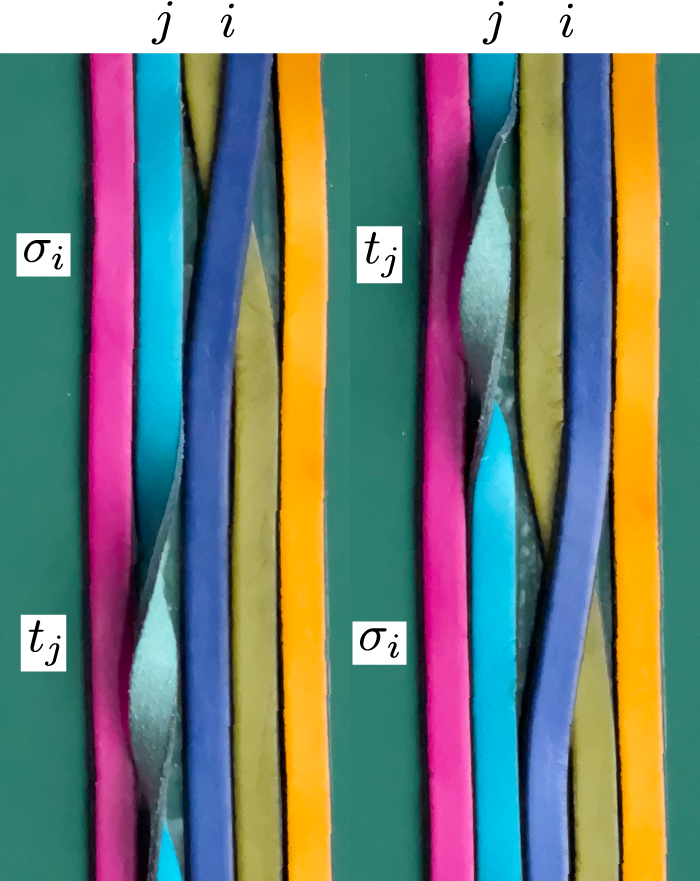}
    \subcaption{$\sigma_i t_j = t_{j}\sigma_i$, \emph{if} $j\notin \{i,i+1\}$}
    \label{fig:fbg_r2}
   \end{minipage}

\vspace{-5pt}
\caption{{Relations for the framed braid group.}}
\label{fig:fbg_relations}
\end{figure}

To describe this mathematically, we turn to the \emph{framed braid groups} $FB_n^{\mathbb{R}^2}$; for more details, we recommend~\cite{ko1992framedbraids}.
If we restrict to framed braids with no twisting in the strands, we obtain an embedding $B_n^{\mathbb{R}^2}\hookrightarrow FB_n^{\mathbb{R}^2}$ where the inclusion of the generators $\sigma_i$ of $B_n$ is  a crossing of strand $i+1$ over strand $i$ with no twists in any strand. The structure of the framed braid group is a semidirect product $\ZZ^n \rtimes B_n^{\RR^2}$, where we take $\{t_i\}_{i=1}^{n}$ as a basis for $\mathbb{Z}^n$ with $t_i$ corresponding to a full twist on the $i\th$ strand, as in Figure~\ref{fig:fbg_gen}. The relations between the framing generators and crossing generators are given in Figure~\ref{fig:fbg_relations}. The \emph{framed pure braid group} $FPB_n^{\RR^2}$ on $n$ strands is defined --- similar to the pure braid group --- as consisting of framed braids where the orders of strands at each end are not permuted with respect to one another.

\subsection*{Flypes} 

A similar kind of object to a magic braid, called a \emph{mexican plait}, was considered by Murasugi and Kurpita in~\cite{murasugi1999}. They gave a list of relations for these mexican plaits, which do not include framing. If we add framing to these relations the result is a different but equivalent set of relations which we call \emph{flype relations}. \linebreak

\vspace{-16pt}
\begin{wrapfigure}[6]{l}{2.25in}
\vspace{-10pt}
    \centering
     \includegraphics[width=2.25in]{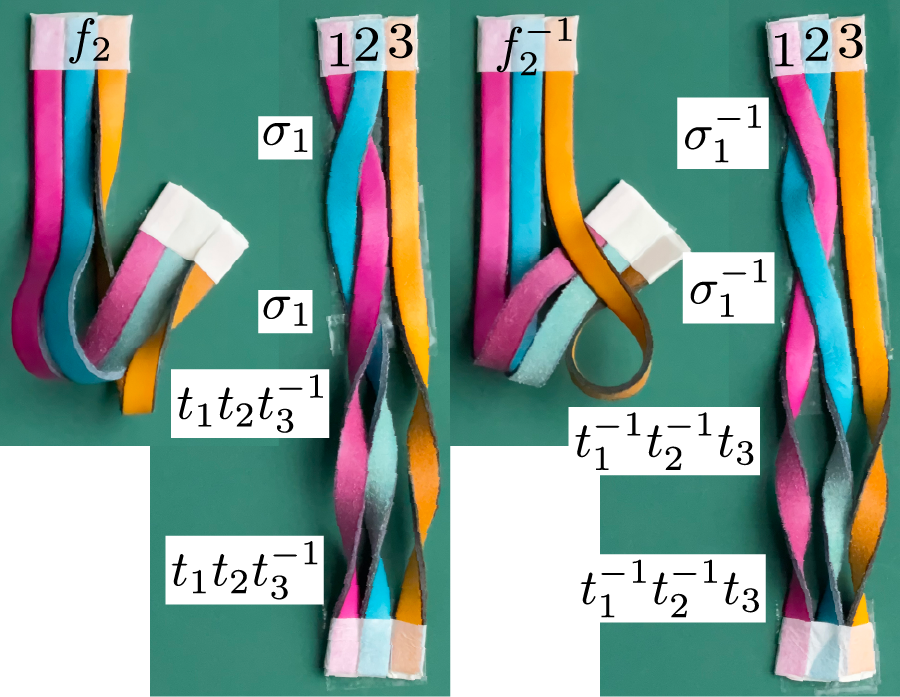}
\vspace{-18pt}
\caption{{Flypes.}}\label{fig:flype}
\end{wrapfigure}\noindent They are as follows:

\vspace{-22pt}
\begin{eqnarray}
 f_1 \hspace{-6pt}&=& \hspace{-6pt} t_1(t_2^{-1} t_3^{-1}\dots t_n^{-1}) (\sigma_{2}^{-1}\sigma_{3}^{-1}\dots\sigma_{n-1}^{n-1})^{n-1},\nonumber\\
f_2  \hspace{-6pt}&=& \hspace{-6pt}  (t_1t_2) (t_3^{-1}\dots t_n^{-1})(\sigma_1)^{2}(\sigma_{3}^{-1}\sigma_{4}^{-1}\dots\sigma_{n-1}^{-1})^{n-2},\nonumber\\
f_3  \hspace{-6pt}&=& \hspace{-6pt}  (t_1t_2 t_3)(t_4^{-1}\dots t_n^{-1})(\sigma_2\sigma_1)^3(\sigma_4^{-1}\sigma_5^{-1}\dots\sigma_{n-1}^{-1})^{n-3},\nonumber\\
\vdots & & \nonumber\\
f_i  \hspace{-6pt}&=& \hspace{-6pt}  (t_1\dots t_i)(t_{i+1}^{-1} \dots t_n^{-1})(\sigma_{i-1}\dots\sigma_{2}\sigma_1)^i(\sigma_{i+1}^{-1}\sigma_{i+2}^{-1}\dots \sigma_{n-1}^{-1})^{n-1},\nonumber\\
\vdots & & \nonumber\\
\hspace{-10pt}f_{n-1}  \hspace{-6pt}&=& \hspace{-6pt}  (t_1\dots t_{n-1})t_n^{-1}(\sigma_{n-2}\sigma_{n-3}\dots\sigma_{2}\sigma_{1})^{n-1}.
\end{eqnarray}

\vspace{-5pt}
\noindent{Figure~\ref{fig:flype} shows two of these flypes. These relations correspond closely with the moves that are, in practice, used to make magic braids: the flype $f_i$ is performed by pulling the bottom end of the braid through the gap between the $i\th$ and $(i+1)\th$ strand from the \emph{back} to the \emph{front}. However, though flypes are key to the creation of magic braids, there is a more efficient way of describing the underlying mathematics.}

\section*{The Spherical Braid Group and Magic Braids}

\noindent If we consider our braids on $n$ strands to lie in the thickened sphere $\mathbb{S}^2\times[0,1]$ rather than the thickened plane, we obtain the \emph{spherical braid group} $B_n^{\mathbb{S}^2}$. The standard presentation of the spherical braid group is the Artin presentation of the planar braid group with one additional relation, 
\begin{equation}
\sigma_1\sigma_2 \dots \sigma_{n-2}\sigma_{n-1}\sigma_{n-1}\sigma_{n-2}\dots \sigma_2\sigma_{1} = \textrm{id},\label{eq:sph_rel}
\end{equation} 
\noindent shown in Figure~\ref{fig:sph_bg} in the three strand case. This relationship between $B_n^{\mathbb{R}^2}$ and $B_n^{S^2}$ is a direct consequence of the fact that the plane embeds in the $2$-sphere.

\begin{figure}[htb!]
\centering
\captionsetup[subfigure]{justification=centering}
\begin{minipage}{2.32 in}
\vspace{-12pt}
\includegraphics[height=1.25in]{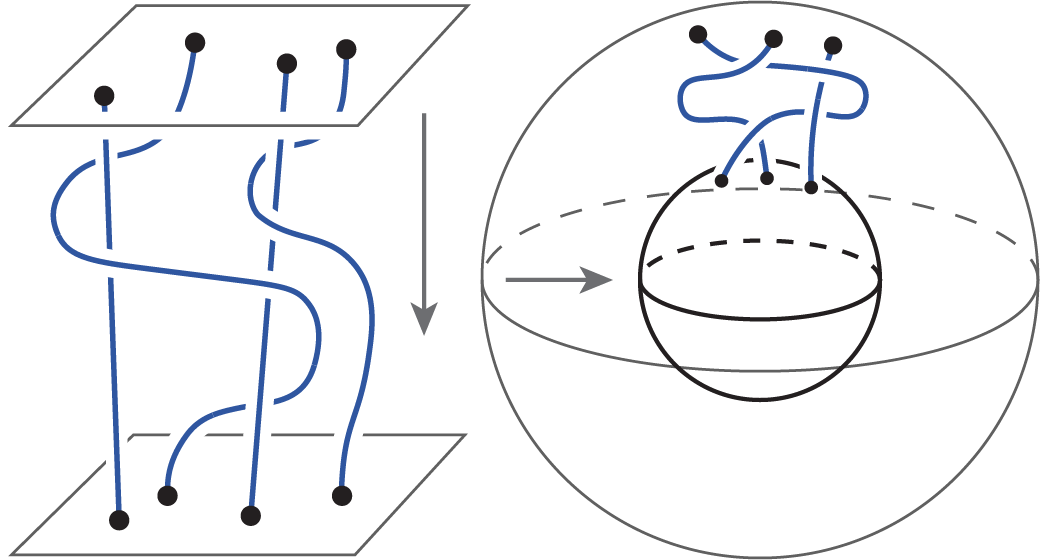}
\subcaption{{Braids whose base space is} $\RR^2$ {vs braids whose base space is} $S^2$.}
\label{fig:r2_vs_s2}
\end{minipage}
\ \
\begin{minipage}{1.75in}
\includegraphics[height=1.25in]{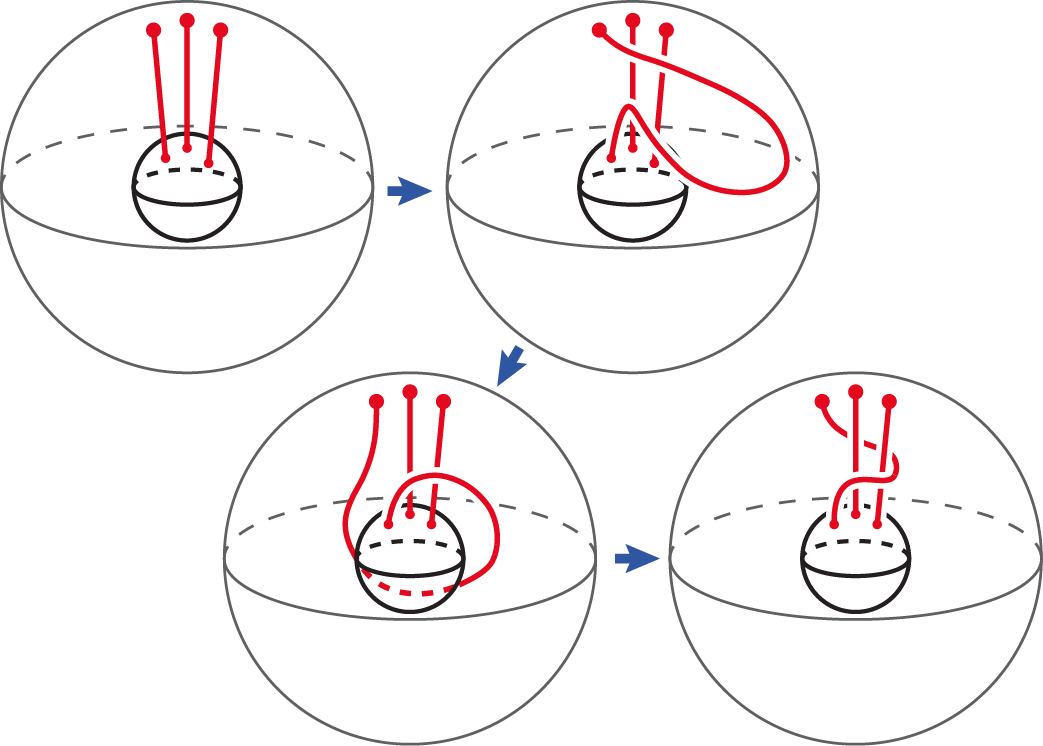}
\subcaption{{The additional relation that distinguishes the spherical braid group from} $B_n$.}
\label{fig:sph_bg}
\end{minipage}
\ \
\begin{minipage}{2.12in}
\vspace{-12pt}
\includegraphics[height=1.25in]{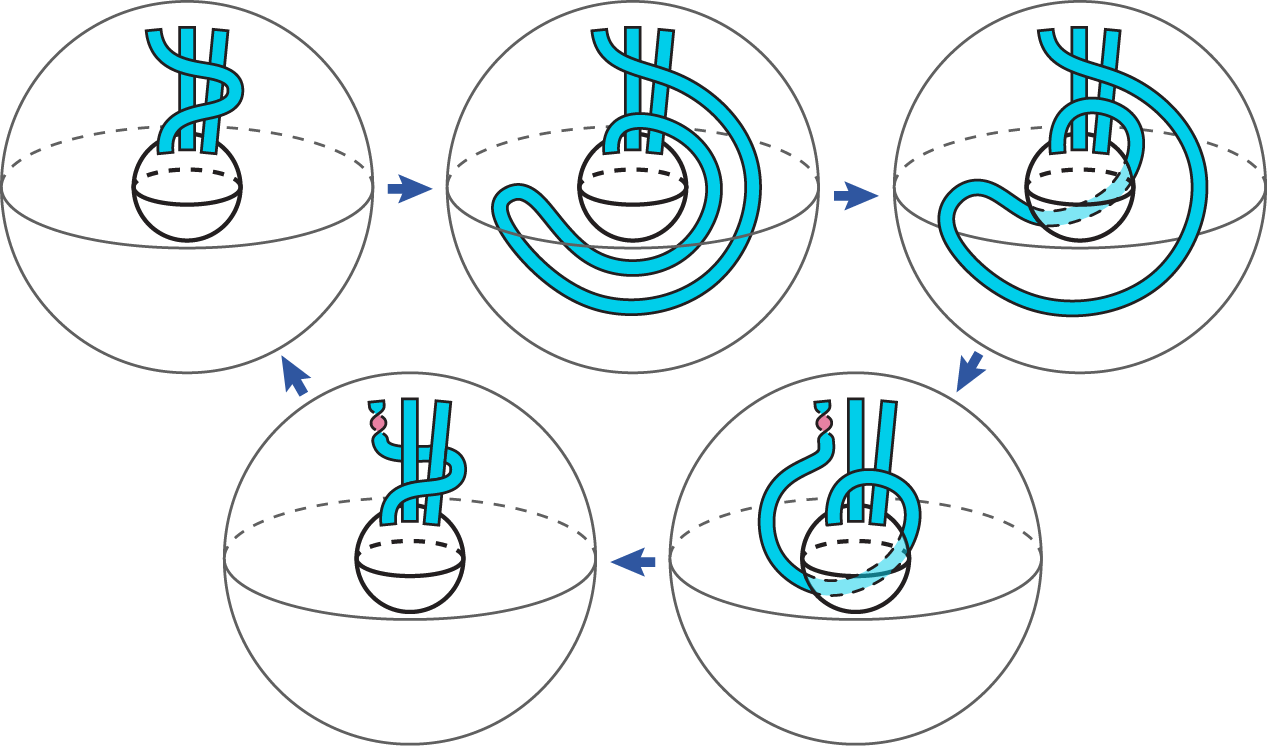}
\subcaption{{The relation shown in Figure 9(b) with the addition of framing.}}
\label{fig:sph_bg}
\end{minipage}

\caption{{Spherical braids.}}
\label{fig:sph}
\end{figure}

We also use framings to describe twisting in spherical braids, as the \emph{spherical framed braid group} $FB_n^{\mathbb{S}^2}$. For more details we recommend~\cite{surfaceframedbraids}.
Similarly to the unframed case, the spherical framed braid group is a quotient of the planar framed braid group by a single additional relation written in Equation~\ref{eq:sph_fr_rel2}, which will also be denoted by $s_{F}$. 
Figure~\ref{fig:mbg_gen}
shows the geometric intuition for why the ribbon embedding this corresponds to is ambient isotopic to the identity framed braid in the thickened sphere. We may also define, similar to the planar case, the \emph{spherical pure braid group} $PB_n^{\mathbb{S}^2}$ and, more pertinently, the \emph{spherical framed pure braid group} $FPB_n^{\mathbb{S}^2}$ of (framed) braids in the thickened sphere $\mathbb{S}^2$ that do not permute their strands from one end to the other.
In the spherical framed pure braid group, the relation $s_F$ can be written in terms of the pure braid group generators $A_{ij}$:

\begin{equation}\label{eq:sph_fr_rel2}
s_F \coloneq t_1^2\sigma_1\sigma_2 \dots \sigma_{n-2}\sigma_{n-1}\sigma_{n-1}\sigma_{n-2}\dots \sigma_2\sigma_{1} =t^2 A_{1 2} A_{1 3} \dots A_{1 n} =  \textrm{id}.
\end{equation} 

\subsection*{The Magic Braid Group}

\begin{wrapfigure}[7]{l}{2in}
\vspace{-100pt}\centering
\includegraphics[height=2in]{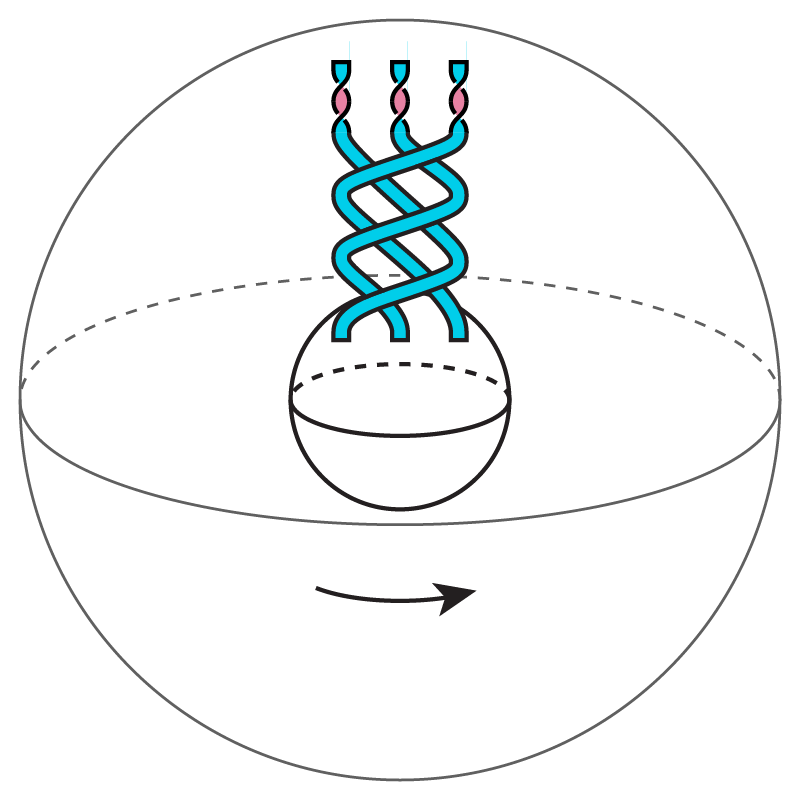}
\caption{{The full twist relation on $3$-strands.}}
\label{fig:fulltwistframed}
\end{wrapfigure}

\noindent The leather objects that we are interested in can be described in terms of framed braid words by taking the twisted ribbons that a framed braid corresponds to and gluing \emph{coupons} to the top and bottom of the braid. Each coupon fixes an order of the ribbons, and any isotopy the leatherworker performs must preserve this order. This gives us the full twist relation, described in Equation~\ref{eq:fulltwistframed} and depicted in Figure~\ref{fig:fulltwistframed}.

We refer to the group given by taking the quotient of $FB_n^{S^2}$ by the full twist relation $\Delta$ as the \emph{leather braid group}, $LB_n$:

\begin{equation}\label{eq:fulltwistframed}
\Delta:= \prod_{i=1}^{n}t_i (\prod_{i=1}^{n-1}\sigma_i)^n.
\end{equation}

To show that quotienting $FB_n^{\mathbb{R}^2}$ by $s_{F}$ and $\Delta$ yields the same group as quotienting by flypes, we give a transposition between $\{f_i\}_{i=1}^{n-1}$ and the relations $s_{F},\Delta$ in Equation~\ref{eq:transposition}. The following transposition shows that any configuration of the cut rectangle of leather that a leatherworker can achieve can be made either by the moves encoded in the relations $s_{F},\Delta$ or by flypes:

 \vspace{-12pt}
\begin{eqnarray}\label{eq:transposition}
  \qquad\qquad\qquad \qquad\qquad\qquad \qquad f_i &=&  \Delta \prod_{j=1}^{i} \left(P_j^{-1}\ldots s_{F} P_j\right),  \textrm{where } P_j:= \prod_{k=1}^{j-1} \sigma_k^{-1} \textrm{ and } P_1=id\nonumber\\
    s_{F} &=& \sigma_1 f_1^{-1} f_2 \sigma_1^{-1} \nonumber\\
    \Delta &=& \sigma_1 f_1^{-1} f_2 \sigma_1^{-1} f_1^{-1} .
\end{eqnarray}

\noindent We encourage the reader to verify these equalities with their own cut leather rectangles.

With this setup we are now able to define the \emph{magic braid group} $MB_n$. Let the map $m:B_n^{\mathbb{R}^2}\rightarrow LB_n$ be the composition of the inclusion $B_n^{\mathbb{R}^2}\hookrightarrow FB_n^{\mathbb{R}^2}$ of the planar braid group into the framed braid group, with the quotients $FB_n^{\mathbb{R}^2}\twoheadrightarrow FB_n^{S^2}\twoheadrightarrow LB_n$ between the framed planar, spherical, and leather braid groups. We define the magic braid group $MB_n$ to be the kernel of $m$. The reason that this is a sensible definition is that the image of $B_n^{\mathbb{R}^2}\hookrightarrow FB_n^{\mathbb{R}^2}$ is the set of ribbon braids without any twists and that the kernel of $FB_n^{\mathbb{R}^2}\rightarrow LB_n$ is the set of framed braids that a leatherworker can manipulate the trivial cut rectangle of leather into. We note that the magic braid group is contained in the spherical framed pure braid group $FPB_n^{\mathbb{S}^2}$, because as noted previously, both ends of a magic braid are fixed in a trivial permutation.

\section*{The Word Problem and New Magic Braids}

\begin{figure}[b]
\centering
\includegraphics[height=2in]{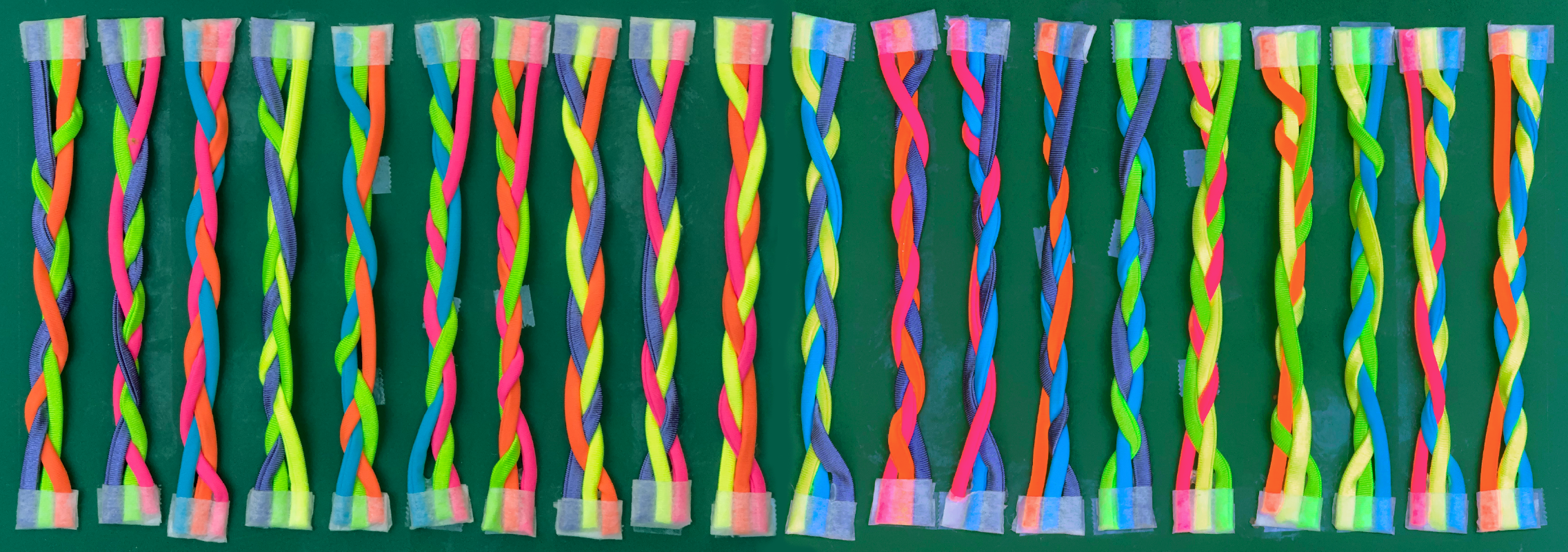}
\caption{{Some three and four strand magic braids discovered using our solution to the word problem of} $LBr_n$.}
\label{fig:newmagicbraids}
\end{figure}
One last observation is that by holding one strand fixed, we can see the spherical framed pure braid group on $n$ strands $FPB_n^{\mathbb{S}^2}$ is a subgroup of the planar framed pure braid group on $n-1$ strands $FPB_{n-1}^{\mathbb{R}^2}$. The planar framed pure braid group on $n-1$ strands splits as $\mathbb{Z}^{n-1}\times PB_{n-1}^{\mathbb{R}^2}$, a product of a group with a trivial word problem with another whose word problem has several well known solutions in the literature. 

We have implemented a program in Python to determine whether leather braids are trivial and search for new magic braids. This led to the discovery of several new magic braids. Figure~\ref{fig:newmagicbraids} shows an arbitrarily chosen selection of three strand braids with around 15 crossings which were fabricated using three different colors for each strand to more clearly show the crossings. The strips were joined with hot glue prior to braiding to keep in the spirit of leather work. These samples benefitted greatly from elastic ribbons, as the effective length of each strand is not consistent. This is because the number of crossings between color pairs can vary within a braid. A larger selection of three- and four-strand examples are presented in the supplementary material. 

\vspace{-5pt}
\section*{An Infinite Family of Magic Braids}

\begin{figure}[b!]
\centering
\captionsetup[subfigure]{justification=centering}
\begin{minipage}{3.5in}
\includegraphics[width=\linewidth]{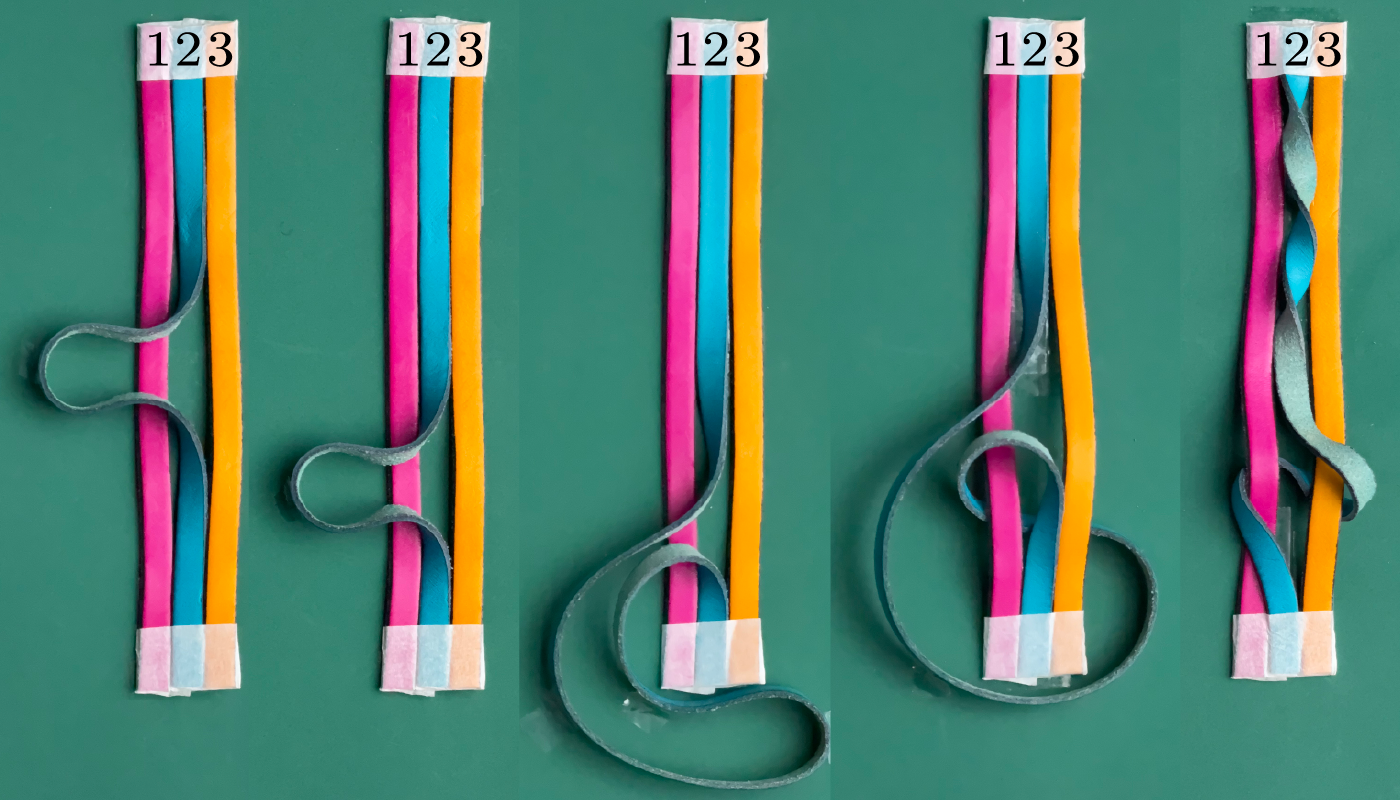}
\subcaption{$h_2$}
\label{fig:mbg_h}
\end{minipage}

\vspace{5pt}
\begin{minipage}{1.1in}
\includegraphics[width=\linewidth]{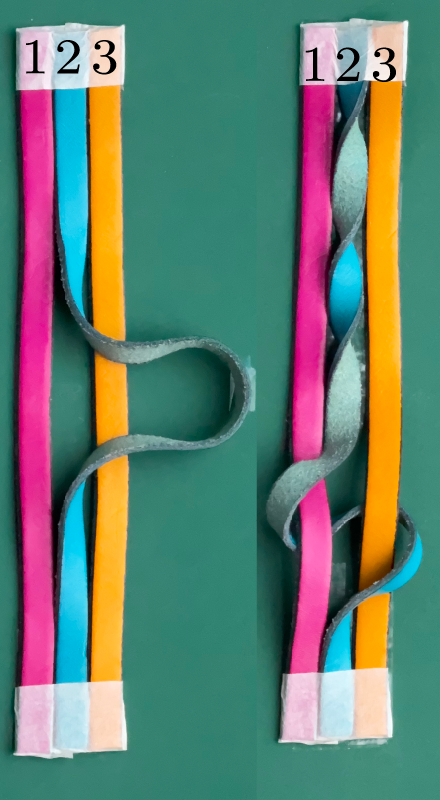}
\subcaption{ $h_2^{-1}$ }
\label{fig:mpg_h_inv}
\end{minipage}
\
\begin{minipage}{1.1in}
\includegraphics[width=\linewidth]{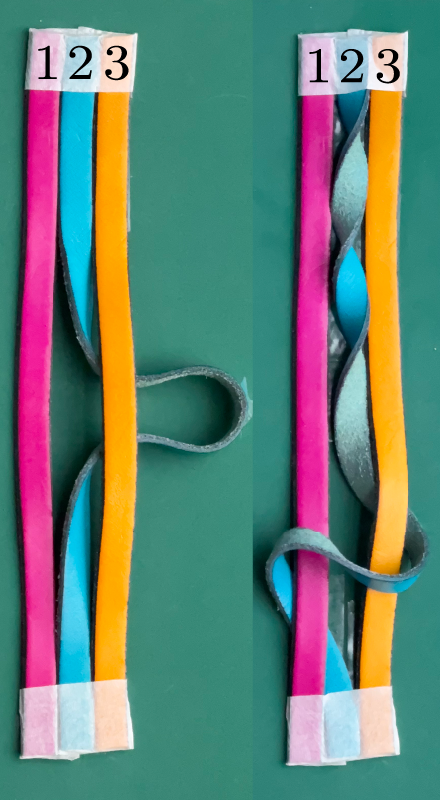}
\subcaption{ $g_2$}
\label{fig:mpg_g}
\end{minipage}
\
\begin{minipage}{1.1in}
\includegraphics[width=\linewidth]{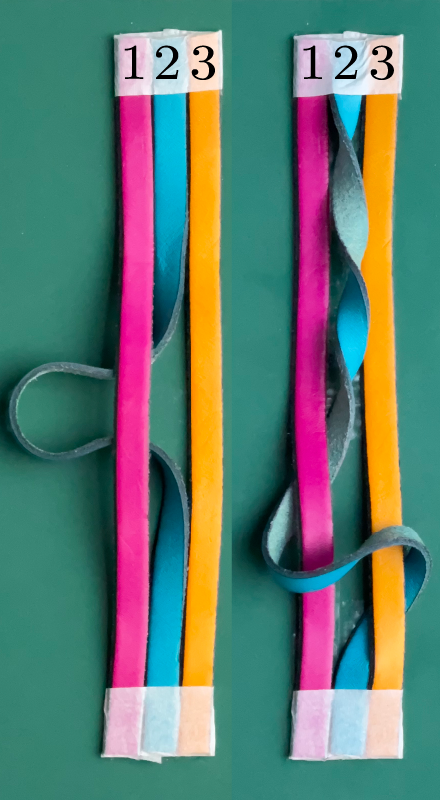}
\subcaption{ $g_2^{-1}$}
\label{fig:mpg_g_inv}
\end{minipage}

\vspace{5pt}
\caption{{Examples of the proposed moves used to generate the magic braid group.}}
\label{fig:mbg_gen}
\vspace{-5pt}
\end{figure}

\begin{figure}[h!]
\centering
\includegraphics[height=2in]{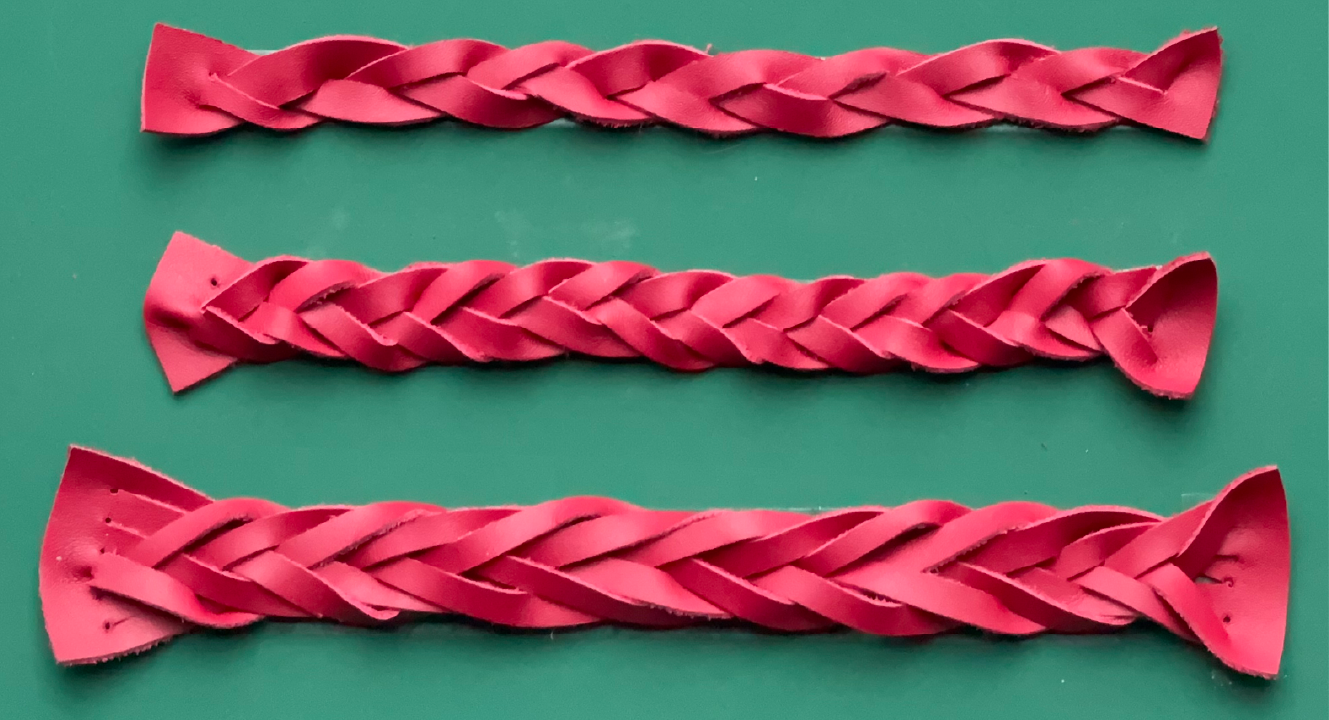}
\caption{{An example of fishtail braids on $3$, $5$, and $7$ strands.} }
\label{fig:fishtails}
\end{figure}

The three strand alternating braid is the start of an infinite family of fishtail braids with one example in each magic braid group with an odd number of strands. Examples of the first three members of this family are shown in Figure~\ref{fig:fishtails}. Note that the full-twist relation was not actually required to form the three strand alternating braid. This turns out to be a general fact about all magic braids. A proof of this and the fact that the fishtail braids are an infinite family of magic braids will be contained in our forthcoming paper.

We will now describe how to make a three strand magic braid by first defining two families of framed braids $g_i$ and $h_i$ for $i\in \{1,\ldots,n\}$. 

\begin{equation}\label{eq:sph_fr_rel}
\begin{array}{ll}
g_i \hspace{-6pt}&\coloneq  t_i^2 A_{i, i+1} A_{i, i+3} \dots A_{i, n} A_{1, i} A_{2, i} \dots A_{i-1, i}\\
& =  Q_i^{-1} s_{F} Q_i\\
h_i \hspace{-6pt}&\coloneq \ t_i^2 A_{n, i} A_{n-1, i} \dots A_{i+1, i} A_{i, i-1} A_{i, i-2} \dots A_{i, 1}\\
& = P_i^{-1} s_{F} P_i\\
\end{array}
\end{equation}

\noindent The first, $g_i$, wraps the $i\th$ ribbon around the bottom of the braid clockwise from the \emph{back}, when viewed from the top of the braid. The second word $h_i$ corresponds to wrapping the $i\th$ ribbon around the bottom of the braid clockwise from the \emph{front}. In both cases the ambient isotopy performed introduces a double twist on the $i$th strand. These words $g_i$ and $h_i$ can be expressed in terms of the crossing and twist generators in the framed planar braid group or as conjugates of the spherical framed braid relation $s_{F}$ as in Equation~\ref{eq:sph_fr_rel}. Recall the definition of $P_j$ given in Equation~\ref{eq:transposition} and let $Q_j := \prod_{i=1}^{j-1}\sigma_j$. These words correspond to the isotopies shown in Figure~\ref{fig:mbg_gen}.

In leatherworking, the most common magic braid is the alternating braid on three ribbons. Mathematically, a set of instructions for how to create this braid is equivalent to an expression of its inclusion in the framed braid group of the plane in terms of the relations that generate the kernel of $FB_n^{\mathbb{R}^2}\rightarrow LB_n$. The construction of the alternating three strand braid can therefore be written as $g_2 h_2^{-1}$ where the words $g_2$ and $h_2$ are themselves conjugates of $s_{F}$. We suspect that all magic braids may be expressed this way as a product of the words $g_i$ and $h_i$.

\vspace{-0pt}
\section*{Conclusion}
Previous work on magic braids has focused on using flypes --- the actual moves used by leatherworkers to create these braids --- to describe the underlying mathematics, but the representation of flypes in braid groups is often complicated. Representing magic braids instead as a kernel in the spherical framed braid group reframes the problem as an isotopy in a thickened sphere, capturing the extra moves and restrictions of magic braids as compared to standard braid groups with simpler relations. This has the potential to allow the easier application of tools from both braid groups and topological theory to investigate the structure of these remarkable objects.

\end{document}